\documentclass{siamltex}
\usepackage{amsmath,amssymb,hyperref,url}

\newcommand{\colspan}{\mathrm{colspan}}
\newcommand{\geqs}{\geqslant}
\newcommand{\Gt}{\tilde{G}}
\newcommand{\gt}{\tilde{g}}
\newcommand{\Kk}{\mathcal{K}}

\newcommand{\Rr}{\mathbb{R}}
\newcommand{\Sigmat}{\tilde{\Sigma}}
\newcommand{\Span}{\mathrm{span}}
\newcommand{\Ut}{\tilde{U}}
\newcommand{\Vt}{\tilde{V}}
\newcommand{\yt}{\tilde{y}}

\title{%
Block Krylov subspace exact time integration\\
of linear ODE systems\\
Part~1: algorithm description%
\thanks{This manuscript is created on September~23, 2011.}
}

\author{Mike A. Botchev\thanks{Department of Applied Mathematics,
University of Twente, P.O.~Box 217, 7500~AE Enschede, the Netherlands,
\url{mbotchev@na-net.ornl.gov}.}}
\date{September 23, 2011}

\begin{document}
\maketitle

\begin{abstract}
We propose a time-exact Krylov-subspace-based method for solving
linear ODE (ordinary differential equation) systems of the 
form $y'=-Ay + g(t)$, where $y(t)$ is the unknown function.
The method consists of two stages.  The first stage is an accurate 
polynomial approximation of the source term $g(t)$, constructed 
with the help of the truncated SVD (singular value decomposition).
The second stage is a special residual-based block Krylov subspace 
method.  

The accuracy of the method is only restricted by the 
accuracy of the polynomial approximation and by the error of the 
block Krylov process.  Since both errors can, in principle, be
made arbitrarily small, this yields, at some costs, a time-exact 
method.  
\end{abstract}



\section{Problem formulation}
Consider initial-value problem (IVP)
\begin{equation}
\label{ivp}
\left\{
\begin{aligned}
&y'   = -Ay + g(t),\\
&y(0) = v,  \qquad\qquad t\in[0,T],
\end{aligned}
\right.  
\end{equation}
where $y(t)$ is the unknown vector function, $y:\Rr\rightarrow\Rr^n$,
and the matrix $A\in\Rr^{n\times n}$, vector function $g:\Rr\rightarrow\Rr^n$, 
and vector $v\in\Rr^n$ are given.

Let $\yt(t)\equiv y(t)-v$ (meaning that $\yt(t)=y(t)-v$ for all $t$).  
Note that the function $\yt(t)$ satisfies IVP 
\begin{equation}
\label{ivp2}
\left\{
\begin{aligned}
&\yt'   = -A\yt + \gt(t),\\
&\yt(0) = 0,  \qquad\qquad t\in[0,T],
\end{aligned}
\right.  
\end{equation}
where $\gt(t)\equiv g(t)-Av$.  We will assume that the IVP~\eqref{ivp}
is brought to the equivalent form~\eqref{ivp2} and, for simplicity,
we omit the tilde sign $\tilde{\cdot}$ in~\eqref{ivp2}.

\section{Polynomial approximation}
We now describe the first stage of the method, the best fit polynomial
approximation of the source term $g(t)$.
Choose $s$ points $0=t_1<t_2<\dots<t_{s-1}<t_s=T$ on the time interval $[0,T]$.
The polynomial approximation is based on the truncated SVD 
(singular value decomposition) of the matrix 
$$
\Gt = 
\begin{bmatrix}
g(t_1) & g(t_2) & \dots & g(t_s)  
\end{bmatrix}\in \Rr^{n\times s},
$$
whose columns are samples $g(t_i)$, $i=1,\dots s$, of the vector
function $g(t)$.
More precisely, let 
\begin{equation}
\label{svd}
\Gt=\Ut\Sigmat\Vt^T, \quad \Sigmat=\text{diag}(\sigma_1,\dots,\sigma_s)\in\Rr^{s\times s},
\quad \sigma_1\geqs\dots\geqs\sigma_s\geqs 0,
\end{equation}
be the thin SVD \cite[Section~2.5.4]{GolVanL}, where 
the matrices $\Ut\in\Rr^{n\times s}$ and $\Vt\in\Rr^{s\times s}$ have
orthonormal columns $u_1$, \dots, $u_s$ and $v_1$, \dots, $v_s$,
respectively.
An approximation to $\Gt$ can be obtained by truncating the SVD
as 
\begin{equation}
\label{svd_trunc}
\Gt=\Ut\Sigmat\Vt^T=\sum_{i=1}^s\sigma_iu_iv_i^T\approx 
                  \sum_{i=1}^m\sigma_iu_iv_i^T=U\Sigma V^T,
\quad m<s,
\end{equation}
where 
$\Sigma\in\Rr^{m\times m}=\text{diag}(\sigma_1,\dots,\sigma_m)$
and the matrices $U\in\Rr^{n\times m}$ and 
$V\in\Rr^{s\times m}$ are formed by the first $m$ columns of
$\Ut$ and $\Vt$, respectively.
Denote the obtained approximate matrix by $G=U\Sigma V^T$.
If follows from~\eqref{svd_trunc} 
that the SVD of $\Gt-G$ is readily available as
$\sum_{i=m+1}^s\sigma_iu_iv_i^T$.  Hence, for the 2-norm and 
Frobenius norm of the error 
$\Gt-G$ holds~\cite[Section~2.5.3]{GolVanL}:
\begin{align*}
\|\Gt-G\|_2 = \sigma_{m+1},\quad
\|\Gt-G\|_F^2 = \sigma_{m+1}^2+\dots+\sigma_s^2.  
\end{align*}
Looking at SVD identity~\eqref{svd} columnwise, we see that 
every sample value  $g(t_i)$ of the function $g(t)$ can be
approximated by a linear combination of the vectors $u_1$,
\dots, $u_m$:
$$
\begin{aligned}
g(t_i) &=
    (\sigma_1v_{i1})u_1 + (\sigma_2v_{i2})u_2 + \dots + (\sigma_sv_{is})u_s
\\
       &\approx
    (\sigma_1v_{i1})u_1 + (\sigma_2v_{i2})u_2 + \dots + (\sigma_sv_{im})u_m,
\end{aligned}
$$
where $v_{ij}$ are the entries of the unitary matrix $V$.
Following the approach of~\cite{botchev_ea09}, we consider the coefficients
of these linear combinations, namely $\sigma_jv_{ij}$,  $j=1,\dots, m$, 
as values of some unknown functions $f_j(t)$ at $t_i$.  These functions
can be easily approximated, at a low cost (typically $m\ll n$) 
and with a very high accuracy, by a polynomial fit~\cite{botchev_ea09}.  
This yields
\begin{equation}
\label{best_fit}  
\begin{aligned}
g(t_i) &\approx f_1(t_i)u_1 + f_2(t_i)u_2 + \dots + f_m(t_i)u_m
\\
       &\approx p_1(t_i)u_1 + p_2(t_i)u_2 + \dots + p_m(t_i)u_m.
\end{aligned}
\end{equation}
For simplicity, we assume that all the best-fit polynomials have
the same the order $r$.  Packing the polynomials $p_j(t)$, 
$j=1,\dots,m$, in one polynomial vector function 
$p(t)=(p_1(t),\dots,p_m(t))^T$, we obtain a polynomial
approximation
\begin{equation}
\label{poly_approx}
g(t)\approx U p(t).  
\end{equation}
There are three sources contributing to the approximation error here.  
First, the quality of the approximation is influenced by the choice of 
the sample points
$t_1$, \dots, $t_s$.  
Second, by the number of terms $m$ in the SVD truncation~\eqref{svd_trunc}
and, finally, by the polynomial best fit in~\eqref{best_fit}.
All these errors can be easily controlled when the approximation
is constructed~\cite{botchev_ea09}, thus giving possibility for an 
adaptive approximation procedure.  With~\eqref{poly_approx},
the original initial-value problem~\eqref{ivp2} takes the form
\begin{equation}
\label{ivp3}
\left\{
\begin{aligned}
&y'   = -Ay + U p(t),\\
&y(0) = 0,  \qquad\qquad t\in[0,T],
\end{aligned}
\right.  
\end{equation}
We now introduce a block Krylov subspace method to solve this problem.

\section{Residual-based block Krylov subspace method}
To construct a Kry\-lov subspace block iterative method for solving~\eqref{ivp3},
we use the exponential residual concept described in~\cite{Botchev2010}.
Choosing the initial guess $y_0(t)$ to be a zero vector function,
we see that the corresponding initial residual is
\begin{equation}
\label{res0}
r_0(t) = -Ay_0(t)-y_0'(t) + U p(t) = U p(t).
\end{equation}
We follow the approach of~\cite{Botchev2010}, where the approximate
solution $y_k(t)$ at Krylov iteration $k$ is obtained as
$$
y_k(t) = y_0(t) + \xi_k(t).
$$
Here the vector function $\xi_k(t)$ is the Krylov subspace 
approximate solution of the correction problem
\begin{equation}
\label{ivp_corr}
\left\{
\begin{aligned}
&\xi'   = -A\xi + r_0(t),\\
&\xi(0) = 0,  \qquad\qquad t\in[0,T],
\end{aligned}
\right.  
\end{equation}
Note that if $\xi_k(t)$ solves~\eqref{ivp_corr} exactly then $y_k(t)$
is the sought-after exact solution of~\eqref{ivp3}.
We solve~\eqref{ivp_corr} by projecting it onto a block Krylov subspace
defined as
$$
\Kk_k(A,U) \equiv \Span\left\{U, AU, A^2U, \dots, A^{k-1}U \right\},
$$
with dimension at most $k\cdot m$.  An orthonormal basis for
this subspace can be generated by the block Arnoldi or Lanczos process described
e.g.\ in~\cite{Henk:book,SaadBook}.  The process produces, after $k$ block
steps, matrices 
$$
V_{[k+1]}=\begin{bmatrix}
        V_1 & V_2 & \dots & V_{k+1}
        \end{bmatrix}
        \in\Rr^{n\times (k+1)m},
\qquad
H_{[k+1,k]}\in\Rr^{(k+1)m\times km}.
$$
Here $V_i\in\Rr^{n\times m}$, $V_1$ is taken to be the matrix $U$ produced by 
the truncated SVD~\eqref{svd_trunc}
and $V_{[k+1]}$ has orthonormal columns spanning the Krylov subspace,
namely, 
$$
\colspan(V_{[k]})=\Kk_k(A,U).
$$
The matrix $H_{[k+1,k]}$ is block upper Hessenberg, with $m\times m$ blocks $H_{ij}$,
$i=1,\dots,k+1$, $j=1,\dots,k$.
The matrices $V_{[k+1]}$ and $H_{[k+1,k]}$ satisfy the block Arnoldi (Lanczos) 
decomposition~\cite{Henk:book,SaadBook}
\begin{equation}
\label{Arnoldi}
AV_{[k]} = V_{[k+1]}H_{[k+1,k]} = V_{[k]}H_{[k,k]} + V_{k+1}H_{k+1,k}E_k^T,
\end{equation}
where $H_{k+1,k}$ is the only nonzero block in the last $k+1$ block row
of $H_{[k+1,k]}$
and $E_k\in\Rr^{n\times k}$ is formed by the last $m$ columns of the 
$km\times km$ identity matrix.

The Krylov subspace solution $\xi_k(t)$ is computed as 
$$
\xi_k(t) = V_{[k]}u(t),
$$ 
where $u(t)$ solves the projected IVP
\begin{equation}
\label{ivp_proj}
\left\{
\begin{aligned}
&u'(t) = -H_{[k,k]} u(t) + V_{[k]}^Tr_0(t),\\
&u(0)  = 0,  \qquad\qquad t\in[0,T].
\end{aligned}
\right.  
\end{equation}
Note that 
$$
V_{[k]}^Tr_0(t) = V_{[k]}^TV_1 p(t) = E_1 p(t),
$$
where $E_1\in\Rr^{km\times m}$ is formed by the first $m$ columns 
of the $km\times km$ identity matrix.
Using~\eqref{res0}, \eqref{Arnoldi} and~\eqref{ivp_proj}, we see that for the exponential 
residual $r_k(t)$ of the solution $y_k(t)$ holds
\begin{equation}
\label{res_k}
\begin{aligned}
r_k(t) &= -Ay_k-y_k'+Up(t) = -Ay_0-y_0'-AV_{[k]}u(t)-V_{[k]}u'(t)+Up(t) =\\
       &= r_0(t)-AV_{[k]}u(t)-V_{[k]}u'(t) =\\
       &= r_0(t)-(V_{[k]}H_{[k,k]} + V_{k+1}H_{k+1,k}E_k^T)u(t)- 
                                  V_{[k]}u'(t) =\\
       &= r_0(t) - V_{[k]}(H_{[k,k]}u(t) + u'(t)) - 
                   V_{k+1}H_{k+1,k}E_k^Tu(t) =\\
       &= r_0(t) - V_kE_1p(t) - V_{k+1}H_{k+1,k}E_k^Tu(t) =
                              - V_{k+1}H_{k+1,k}E_k^Tu(t).  
\end{aligned}
\end{equation}
A similar expression for the exponential residual is obtained in~\cite{Botchev2010}
for a non-block Krylov subspace method.  There are two important messages
relation~\eqref{res_k} provides.  First, the residual can be computed efficiently
during the iteration process because the matrices $V_{k+1}$ and $H_{k+1,k}$ 
are readily available in the Arnoldi or Lanczos process.  Second,
the residual after $k$ block steps has the same form as the initial
residual~\eqref{res0}, namely it is a matrix of $m$ orthonormal columns times
a time dependent vector function.  This allows for a restart in the block
Krylov method: set $y_0(t):=y_k(t)$, then relation~\eqref{res0} holds with
$U:=V_{k+1}$ and $p(t):=-H_{k+1,k}E_k^Tu(t)$.  The just described correction
with $k$ block Krylov iterations can then be repeated,
which results in a restarted block Krylov subspace method for 
solving~\eqref{ivp3}.

\bibliography{my_bib}
\bibliographystyle{abbrv}

\end{document}